\documentclass[11pt]{article}
\usepackage{amsmath}
\usepackage{amssymb}
\usepackage{bbm}
\usepackage[mathscr]{euscript}
\usepackage[margin=2cm]{geometry}
\usepackage{hyperref}
\usepackage{tikz}
\usetikzlibrary{decorations.pathmorphing,arrows}
\usepackage{braket}
\usepackage{setspace}

\newtheorem{thm}{Theorem}[section]

\newtheorem{prop}[thm]{Proposition}

\newtheorem{lemma}[thm]{Lemma}

\newtheorem{theorem}[thm]{Theorem}
\newtheorem{remark}[thm]{Remark}
\newtheorem{proposition}[thm]{Proposition}

\newtheorem{corollary}[thm]{Corollary}

\newtheorem{conj}[thm]{Conjecture}

\newenvironment{proofof}{}{\hfill$\square$\vskip.5cm}
\newcommand{\R}{\mathbb{R}}
\newcommand{\N}{\mathbb{N}}

\newcommand{\Z}{\mathbb{Z}

}

\renewcommand{\Z}{{\mathbb{Z}}}

\title{Multifold Convolutions, Generating Functions and 1d Random Walks}
\author{Timothy Li\\
\small Vestavia Hills High School,
2235 Lime Rock Rd, Vestavia Hills, AL 35216\\[10pt]
Shannon Starr\\
\small
University of Alabama at Birmingham,
1402 Tenth Avenue South,
 Birmingham, AL 35294-1241
}
\date{}

\begin{document}

\maketitle

\abstract{We consider multifold convolutions of a combinatorial sequence $(a_n)_{n=0}^{\infty}$:
namely,  for each $k \in  \N$ the $k$-fold convolution is
$\mathcal{M}^{(k)}_n(\boldsymbol{a}) = \sum_{j_1+\dots+j_k=n} a_{j_1} \cdots a_{j_k}$.
Let $C_n$ be the Catalan numbers, and let $B_n$ be the central binomial coefficients.
Then for random Dyck paths or simple random walk bridges,
the multifold convolutions give moments of returns to the origin, using the stars-and-bars
problem.
There are well-known explicit formulas for the multifold convolutions of $C_n$
and $B_n$.
But even for combinatorial sequences $B_n^2$ and $B_n^3$, one may determine asymptotics
of multifold convolutions for large $n$.
We also discuss  large deviations: In a second part of the paper we consider an elementary
version of the circle method for calculating asymptotics using complex analysis.
}

\section{Introduction and first results}

Suppose that one has a combinatorial sequence $\boldsymbol{a} = (a_n)_{n=0}^{\infty}$.
Then let us denote its $k$-fold convolution as another sequence $\boldsymbol{\mathcal{M}}^{(k)}(\boldsymbol{a})$
where
$$
	\mathcal{M}^{(k)}_n(\boldsymbol{a})\, =\, \sum_{j_1=0}^{\infty} \cdots \sum_{j_k=0}^{\infty}
\mathbf{1}_{\{n\}}(j_1+\dots+j_k) \cdot a_{j_1} \cdots a_{j_k}\, .
$$
Then we consider the following set-up. Our intuitive motivation is coming from the field of regenerative
processes.
In particular, an inspiring paper for us is \cite{BasuBhatnagar}.
But in that paper by Basu and Bhatnagar they are interested in establishing a CLT.
For us, we are more interested in explicit simple calculations that come from
the multifold convolution formula.
Therefore, consider the following set-up.

Let $S$ be some set, with a distinguished element $\o$. Suppose that
$a_n = |\mathcal{X}_n|$ for a sequence of sets $\mathcal{X}_n \subset S^n$
satisfying the following conditions:
\begin{itemize}
\item[\bf (i)] $\forall n \in \N = \{1,2,\dots\}$,
$\forall \vec{x} = (x_1,\dots,x_n) \in \mathcal{X}_n$ we have $x_n = \o$;
\item[\bf (ii)] $\forall m,n \in \N$, $\forall \vec{x} = (x_1,\dots,x_m) \in \mathcal{X}_m$,
$\forall \vec{y} = (y_1,\dots,y_n) \in \mathcal{X}_n$, if we define
$\vec{z} \in S^{m+n}$ as $\vec{z} = (x_1,\dots,x_m,y_1,\dots,y_n)$ then $\vec{z}$ is in $\mathcal{X}_{m+n}$;
\item[\bf (iii)] $\forall n \in \N$, $\forall \vec{x} = (x_1,\dots,x_n) \in \mathcal{X}_n$, if there is some
$k<n$ such that $x_k = \o$, then $(x_1,\dots,x_k) \in \mathcal{X}_k$
and $(x_{k+1},\dots,x_n) \in \mathcal{X}_{n-k}$.
\end{itemize}
We assume $\mathcal{X}_0$ is a set of cardinality 1, so that $a_0=1$.
Then for each $\vec{x} \in \mathcal{X}_n$ we define the set
$$
\mathfrak{Z}_n(\vec{x})\, =\, \{0\} \cup \{k \in [n]\, :\,
x_k = \o\}\, .
$$
Then the following is true.
\begin{lemma}
For each $k \in \{0,1,\dots\}$ and $n \in \N$,
$$
	\mathcal{M}^{(k+1)}_n(\boldsymbol{a})\, =\, \sum_{\vec{x} \in \mathcal{X}_n} \binom{|\mathfrak{Z}_n(\vec{x})|+k-1}{k}\, .
$$
\end{lemma}
This follows as in the ``stars-and-bars'' problem of Feller. See Section \ref{sec:proofs}.
Let us consider examples.
\begin{proposition}[Catalan, see also Larcombe and French \cite{LarcombeFrench}]
\label{prop:Catalan}
For each $n \in \{0,1,\dots\}$ let $C_n$ be the $n$th Catalan number 
$$
	C_n\, =\, \binom{2n}{n}\, \frac{1}{n+1}\, .
$$
Then
$$
	\mathcal{M}_n^{(k)}(\boldsymbol{C})\, =\, \frac{k}{2n+k} \binom{2n+k}{n}\, .
$$
\end{proposition}
\begin{corollary}
\label{cor:Dyck}
For each $n \in \N$, let $\mathbf{P}^{\mathrm{Dyck}}_n$ denote the uniform measure on the set of all Dyck paths
of length $2n$: a random walk path $(x_1,\dots,x_{2n})$ such that $x_{2n}=0$
and $x_1,\dots,x_{2n-1}\geq 0$. Then, letting
$$
	\mathfrak{Z_n}\, =\, \{0\} \cup \{t \in [n]\, :\, x_{2t}=0\}\, ,
$$
we have 
$$
	\lim_{n \to \infty} \mathbf{P}_n^{\mathrm{Dyck}}(|\mathfrak{Z}_n|=k)\, =\, (k-1)2^{-k}\, ,
$$
for each $k \in \{2,3,\dots\}$.
\end{corollary}
The distribution above is the sum of two IID Geometric-$1/2$ random variables.
It is also called a negative binomial $(2,p)$ random variable with $p=1/2$.

The following proposition is trivial using generating functions.
\begin{proposition}
\label{prop:central}
For each $n \in \{0,1,\dots\}$ let $B_n$ be the $n$th central binomial coefficient 
$$
	B_n\, =\, \binom{2n}{n}\, .
$$
Then
$$
	\mathcal{M}_n^{(k)}(\boldsymbol{B})\, =\, 4^n\, \cdot \frac{(n-1+(k/2))_n}{n!}\, ,
$$
where $(x)_n$ is the falling factorial Pochhammer symbol $x(x-1)\cdots(x-n+1)$.
\end{proposition}

\begin{corollary}
\label{cor:SRWB}
For each $n \in \N$, let $\mathbf{P}^{\mathrm{SRWB}}_n$ denote the uniform measure on the set of all simple random
walk bridges of length $2n$: all simple random walk paths $(x_1,\dots,x_{2n})$ such that $x_{2n}=x_0=0$.
Then, letting
$$
	\mathfrak{Z}_n\, =\, \{0\} \cup \{t \in [n]\, :\, x_{2t}=0\}\, ,
$$
we have
$$
	\lim_{n \to \infty} \mathbf{E}_n^{\mathrm{SRWB}}\left[\left(n^{-1/2}|\mathfrak{Z}_n|\right)^k\right]\, 
=\, 2^k\, \Gamma\left(\frac{k}{2}+1\right)\, 
$$
for each $k \in \{2,3,\dots\}$.
\end{corollary}
\begin{remark}
The limits of the moments equal
$\mathbf{E}\left[\left(\sqrt{2(\mathsf{X}^2+\mathsf{Y}^2)}\right)^k\right]$,
where $\mathsf{X},\mathsf{Y}$ are IID $\mathcal{N}(0,1)$ random variables.
Since they satisfy Carleman's criterion, this shows the weak limit.
That limit is  well known. See for instance Pitman \cite{Pitman}, equation (3).
\end{remark}
\begin{theorem}
\label{thm:2SRWB}
For the example of a combinatorial sequence
$$
	a_n\, =\, \binom{2n}{n}^2\, ,
$$
we have 
$$
	\mathcal{M}_n^{(k)}(\boldsymbol{a})\, \sim\, \frac{e^{4n \ln(2)}}{n\pi} 
\left(\frac{\ln(n)}{\pi}\right)^{k-1} k\, ,\ \text{ as  $n \to \infty$.}
$$
Hence, letting $\mathbf{P}_n^{\mathrm{2SRWB}}$ denote the IID measure on
two simple random walk bridges of length $2n$, $(x_1^{(i)},\dots,x_{2n}^{(i)})$ for $i=1,2$, and letting
$$
	\mathfrak{Z}_n\, =\, \{0\} \cup \{t \in [n]\, :\, x_{2t}^{(1)}=x_{2t}^{(2)}=0\}\, ,
$$
we have
$$
	\lim_{n \to \infty} \mathbf{E}_n^{\mathrm{2SRWB}}\left[\left(\frac{\pi}{\ln(n)}\, \mathfrak{Z}_n\right)^k\right]\, 
	=\, (k+1)!\, .
$$
\end{theorem}
\begin{remark}
The moments of $\pi \mathfrak{Z_n}/\ln(n)$ converge to the moments of a random variable
$\mathsf{X}$ which is Gamma-(2,1) distributed. Again, Carleman's criterion guarantees
uniqueness of the distribution. Hence, $\pi \mathfrak{Z}_n/\ln(n)$ converges weakly to a 
Gamma-(2,1) random variable.
\end{remark}
\begin{theorem}
\label{thm:3SRWB}
For the example of a combinatorial sequence
$$
	a_n\, =\, \binom{2n}{n}^3\, ,
$$
we have 
$$
	\mathcal{M}_n^{(k)}(\boldsymbol{a})\, \sim\, \frac{e^{6n\ln(2)}}{(n\pi)^{3/2}} 
\left(\frac{\pi}{\Gamma\left(\frac{3}{4}\right)^4}\right)^{k-1} k\, .
$$
Hence, letting $\mathbf{P}_n^{\mathrm{3SRWB}}$ denote the IID measure on
three simple random walk bridges of length $2n$, $(x_1^{(i)},\dots,x_{2n}^{(i)})$ for $i=1,2,3$, and letting
$$
	\mathfrak{Z}_n\, =\, \{0\} \cup \{t \in [n]\, :\, x_{2t}^{(1)}=x_{2t}^{(2)}=x_{2t}^{(3)}=0\}\, ,
$$
we have, letting $p =\Gamma(3/4)^4/\pi$,
$$
	\lim_{n \to \infty} 
\mathbf{P}_n^{\mathrm{3SRWB}}(|\mathfrak{Z}_n|=k)\, =\, (k-1)(1-p)^{k-2}p^2\, ,
$$
for each $k \in \{2,3,\dots\}$.
\end{theorem}
\begin{remark}
The last example shares in common with the first example, that the combinatorial sequence
is asymptotically $Ca^n/n^{3/2}$.
The other examples are of the form $Ca^n/n^p$ for some $p\leq 1$, hence not summable.
For the summable examples, the returns to $\o$ are localized near the left and right endpoints.
That is why one obtains a negative-binomial distribution $(n,p)$ with $n=2$.
\end{remark}

\section{A weak LDP for $\mathfrak{Z}_n$ relative to Dyck and SRWB paths}

We will not aim for the usual level of generality. Instead, we merely note the following:
\begin{theorem}
\label{thm:rate}
\label{thm:LDP}
We have, for any $x \in [0,1]$,
$$
	\lim_{n \to \infty} n^{-1} \ln \left(\mathbf{P}^{\mathrm{SRWB}}_n(\mathfrak{Z}_n \geq nx )\right)\, =\, -I^{\mathrm{SRWB}}(x)\, ,
$$
and
$$
	\lim_{n \to \infty} n^{-1} \ln \left(\mathbf{P}^{\mathrm{Dyck}}_n(\mathfrak{Z}_n \geq nx )\right)\, =\, -I^{\mathrm{Dyck}}(x)\, ,
$$
for
$$
	I^{\mathrm{SRWB}}(x)\, =\, (2-x) \ln(2) + (1-x) \ln(1-x) - (2-x) \ln(2-x)\, ,
$$
and
$$
	I^{\mathrm{Dyck}}(x)\, =\, 2 \ln(2) + (1-x) \ln(1-x) - (2-x) \ln(2-x)\, .
$$
\end{theorem}

\begin{figure}
\begin{center}
\begin{tikzpicture}
	\draw (0,0) node[] {\includegraphics[height=5cm]{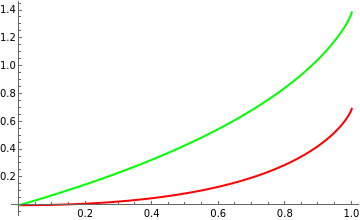}};
\end{tikzpicture}
	\caption{We plot both large deviation rate functions together: red is for the 2SRWB and
green is for the Dyck path. Note that $\ln(2) \approx 0.693147$ and $2\ln(2) \approx 1.38629$.
\label{fig:PREmultipleRate}}
\end{center}
\end{figure}
We show plots of the two rate functions in Figure \ref{fig:PREmultipleRate}.
This theorem will be proved in Section \ref{sec:rates}.
But the proof relies on several calculus-formula lemmas that we state now.

Consider the function
$U : [0,\infty) \times [0,\infty)$ defined as
$$
U(x,y)\, =\, (x+y) \ln(x+y) - x \ln(x) - y \ln(y)\, .
$$
Note that this is an increasing function of $x$ for each $y>0$. 
\begin{lemma}
\label{lem:LDP1}
For each $y \in [0,\infty)$, we have
$$
\lim_{n \to \infty} n^{-1} \ln\left(\mathbf{E}^{\mathrm{SRWB}}_n\left[\exp\left(-n U\left(n^{-1}\mathfrak{Z}_n,y\right)\right)\right]\right)\,
	=\, \mathcal{U}^{\mathrm{SRWB}}(y)\,  ,
$$
and
$$
\lim_{n \to \infty} n^{-1} \ln\left(\mathbf{E}^{\mathrm{Dyck}}_n\left[\exp\left(-n U\left(n^{-1}\mathfrak{Z}_n,y\right)\right)\right]\right)\,
	=\, \mathcal{U}^{\mathrm{Dyck}}(y)\, ,
$$
where
$$
	\mathcal{U}^{\mathrm{SRWB}}(y)\, \stackrel{\mathrm{def}}{:=}\, 
	\left(1+\frac{y}{2}\right)\, \ln\left(1+\frac{y}{2}\right) - \frac{y}{2}\, \ln\left(\frac{y}{2}\right)\, ,
$$
and
$$
	\mathcal{U}^{\mathrm{SRWB}}(y)\, \stackrel{\mathrm{def}}{:=}\, 
\left(2+y\right)\, \ln\left(2+y\right) - (1+y)\, \ln\left(1+y\right) - 2 \ln(2)\, .
$$
\end{lemma}

Since $U(\cdot,y)$ is increasing for a fixed $y$, it may be used in combination with Chebyshev's inequality.
Let us gather some elementary calculations in a lemma.
\begin{lemma}
\label{lem:CalFacts1}
For each $x \in [0,1]$ define the functions $V^{\mathrm{SRWB}}(x,\cdot) : [0,\infty) \to \R$ 
and $V^{\mathrm{Dyck}}(x,\cdot) : [0,\infty) \to \R$ as 
$$
	V^{\mathrm{SRWB}}(x,y)\, =\, \mathcal{U}^{\mathrm{SRWB}}(y) - U(x,y)\ \text{ and }\
	V^{\mathrm{Dyck}}(x,y)\, =\, \mathcal{U}^{\mathrm{Dyck}}(y) - U(x,y)\, .
$$
Then for $x \in (0,1)$ fixed, both functions possess unique minimizers, located at
$$
Y_*^{\mathrm{SRWB}}(x)\, =\, \frac{x^2}{2-2x}\ \text{ and }\ 
Y_*^{\mathrm{Dyck}}\, =\, \frac{x}{1-x}\, ,
$$ 
respectively. Moreover, 
\begin{equation}
\label{eq:IRateSub}
	V^{\mathrm{SRWB}}(x,Y_*^{\mathrm{SRWB}}(x))\, =\, -I^{\mathrm{SRWB}}(x)\ \text{ and }\
	V^{\mathrm{Dyck}}(x,Y_*^{\mathrm{Dyck}}(x))\, =\, -I^{\mathrm{Dyck}}(x)\, .
\end{equation}
\end{lemma}
Thus, using Chebyshev's inequality and these calculus facts we may prove the upper bounds
$$
	\limsup_{n \to \infty} n^{-1} \ln \left(\mathbf{P}^{\mathrm{SRWB}}_n(\mathfrak{Z}_n \geq nx )\right)\, \leq\, -I^{\mathrm{SRWB}}(x)\, ,
$$
and
$$
	\limsup_{n \to \infty} n^{-1} \ln \left(\mathbf{P}^{\mathrm{Dyck}}_n(\mathfrak{Z}_n \geq nx )\right)\, \leq\, -I^{\mathrm{Dyck}}(x)\, .
$$
But, defining the functions $W^{\mathrm{SRWB}}(\cdot,y) : [0,1] \to \R$ and $W^{\mathrm{Dyck}}(\cdot,y) : [0,1] \to \R$ via
$$
W^{\mathrm{SRWB}}(\cdot,y)\, =\, U(x,y) - I^{\mathrm{SRWB}}(x)\ \text{ and }\ 
W^{\mathrm{Dyck}}(\cdot,y)\, =\, U(x,y) - I^{\mathrm{Dyck}}(x)\, ,
$$
we have the following calculus style formulas.
\begin{lemma}
\label{lem:CalFacts2}
For each $y >0$ fixed, the functions $W^{\mathrm{SRWB}}(\cdot,y)$ and $W^{\mathrm{Dyck}}(\cdot,y)$ are strictly concave in $x$ (for $x \in [0,1]$).
Moreover, defining the numbers $X_*^{\mathrm{SWRB}}(y)$ and $X_*^{\mathrm{Dyck}}(y)$ by
$$
	X_*^{\mathrm{SRWB}}(y)\, =\, \sqrt{y(y+2)} - y\ \text{ and }\ 
	X_*^{\mathrm{Dyck}}(y)\, =\, \frac{y}{1+y}\, ,
$$ 
we have that these are the maximizers of the two functions, respectively. (The derivatives are zero.)
Moreover, 
\begin{equation}
\label{eq:Wsub}
W^{\mathrm{SRWB}}(X_*^{\mathrm{SWRB}}(y),y)\, =\, \mathcal{U}^{\mathrm{SWRB}}(y)\ \text{ and }\ 
W^{\mathrm{Dyck}}(X_*^{\mathrm{Dyck}}(y),y)\, =\, \mathcal{U}^{\mathrm{Dyck}}(y)\, .
\end{equation}
\end{lemma}

This is similar for the usual Cram\'er's rule of noting that the Legendre transform of the Legendre transform of the moment generating function
returns the original function.
With all of these calculus formula established, the proof of Theorem \ref{thm:LDP} will be stated
in Section \ref{sec:rates}.
The proofs of the lemmas themselves will be presented in Section \ref{sec:LEMrates}. 
All that is necessary is the lower bound.
By the calculus formulas, if the lower bound is violated at any point then this would lead to a contradiction of Lemma \ref{lem:LDP1} at the corresponding point.
In all of this, we follow the well-known general method of large deviations theory.
For example, we refer to Chapter 2 of the excellent notes of Rezakhanlou \cite{Rezakhanlou}.

\subsection{Some observations about the LDP}

Firstly, note that the two large deviation rate functions differ just as
$$
I^{\mathrm{Dyck}}(x) - I^{\mathrm{SRWB}}(x)\, =\, x\, \ln(2)
$$
which aligns with intuition coming from
the excursion construction of these types of paths.

Let us note the consistency at the endpoints.
The right-derivative of $I^{\mathrm{Dyck}}(x)$ at $x=0$ is $\ln(2)$.
By Corollary \ref{cor:Dyck}, we can see that 
$$
\lim_{n \to \infty} \mathbf{P}^{\mathrm{Dyck}}(|\mathfrak{Z}_n|\geq k)\, \sim\, (k-1)2^{-(k-1)}\, ,
$$
as $k \to \infty$. Thus 
$$
\lim_{k \to \infty} \lim_{n \to \infty} \frac{1}{k} \ln\left(\mathbf{P}^{\mathrm{Dyck}}(|\mathfrak{Z}_n|\geq k)\right)\, 
=\, - \ln(2)\, ,
$$
which is consistent.

For the SRWB, Corollary \ref{cor:SRWB} implies that the $n \to \infty$ scaling limit
of $n^{-1/2} |\mathfrak{Z}_n|$ has cdf which is a Gaussian pdf modulo a constant normalization.
The logarithm of the tail probability behaves like $-R^2/4$ (for the probability to be greater than $R$).
If we take $R \sim \rho n^{1/2}$ we would get $-\rho^2 n/4$.
And this is consistent with $I^{\mathrm{SRWB}}(x) \sim -x^2/4$ as $x \to 0^+$.

We also note for later reference that 
\begin{equation}
	\mathbf{P}^{\mathrm{Dyck}}(\mathfrak{Z}_n=n)\, =\, \frac{1}{C_n}\ \text{ so }
	\lim_{n \to \infty} \frac{1}{n}\, \ln\left(\mathbf{P}^{\mathrm{Dyck}}(\mathfrak{Z}_n=n)\right)\,
	=\, -2\ln(2)\, .
\end{equation}
That is because the only path $\vec{x} \in \Z^{2n}$ that has $\mathfrak{Z}_n(\vec{x})=n$ is $(1,0,1,0,\dots,1,0)$.

In comparison, for the 1SRWB model, the paths that have $\mathfrak{Z}_n(\vec{x})$
are any $(\sigma_1,0,\sigma_2,0,\dots,\sigma_n,0)$ for $\sigma_1,\dots,\sigma_n \in \{1,-1\}$.
So
\begin{equation}
	\mathbf{P}^{\mathrm{SRWB}}(\mathfrak{Z}_n=n)\, =\, \frac{2^n}{B_n}\ \text{ so }
	\lim_{n \to \infty} \frac{1}{n}\, \ln\left(\mathbf{P}^{\mathrm{SRWB}}(\mathfrak{Z}_n=n)\right)\,
	=\, -\ln(2)\, .
\end{equation}
In Figure \ref{fig:PREmultipleRate} we can see these observations reflected
in the right-endpoint values.

We will give an implicit formula for the large deviation rate function of the 2SRWB model, by
a more involved argument.
Therefore, let us also note that
\begin{equation}
\label{eq:right2SRWB}
	\mathbf{P}_2^{\mathrm{SRWB}}(\mathfrak{Z}_n=n)\, =\, \frac{2^{2n}}{B_n^2}\ \text{ so }
	\lim_{n \to \infty} \frac{1}{n}\, \ln\left(\mathbf{P}_2^{\mathrm{SRWB}}(\mathfrak{Z}_n=n)\right)\,
	=\, -2\ln(2)\, .
\end{equation}
This is because, if we write $\vec{x} = (x_1,x_2,\dots,x_{2n})$ for $x_1,x_2,\dots,x_{2n} \in \Z^2$,
then the paths that contribute to the event $\mathfrak{Z}_n(\vec{x})=n$ are of the form
$((\sigma_1^{(1)},\sigma_1^{(2)}),(0,0),\dots,(\sigma_n^{(1)},\sigma_n^{(2)}),(0,0))$
for 
$$
\sigma_1^{(1)},\sigma_1^{(2)},\dots,\sigma_n^{(1)},\sigma_n^{(2)} \in \{1,-1\}\, .
$$

\section{Implicit formulation of the rate function for the 2SRWB}

In this section, we will simply assume that there is a large deviation principle, with an unknown rate function.
Then we will identify the rate function. Hence, let us make the following ansatz.
\begin{conj}
We have, for each $x \in [0,1]$,
$$
	\lim_{n \to \infty} n^{-1} \ln \left(\mathbf{P}^{\mathrm{2SRWB}}_n(\mathfrak{Z}_n \geq nx )\right)\, =\, -I^{\mathrm{2SRWB}}(x)\, ,
$$
for a regular function to be determined.
\end{conj}

Let us denote the complete elliptic integrals of the first and second kind as 
$$
	\mathscr{K}(z)\, =\, \int_0^1 \frac{1}{\sqrt{(1-x^2)(1-zx^2)}}\, dx\, .
$$
and
$$
	\mathscr{E}(z)\, =\, \int_0^1 \frac{\sqrt{1-zx^2}}{\sqrt{1-x^2}}\, dx\, ,
$$
for $z \in [0,1]$.
These are not the most conventional definitions. They are the conventions that are used by Wolfram Mathematica, which is convenient for 
obtaining calculus-style formulas and visualizing plots.
With this, we may define the rate function implicitly.
\begin{prop}
\label{prop:2SRWBrate}
For $t \in [0,1/16]$ define
$$
	\widetilde{X}(t)\, =\, \frac{2(1-16t) \mathscr{K}(16t)}{\mathscr{E}(16t)-(1-16t)\mathscr{K}(16t)}\, \left(\frac{2}{\pi}\, \mathscr{K}(16t)-1\right)\, .
$$
Then we have 
$$
I^{\mathrm{2SRWB}}\left(\widetilde{X}(t)\right)\, =\, \widetilde{\mathcal{I}}(t)\, ,
$$
for the function
$$
\widetilde{\mathcal{I}}(t)\, \stackrel{\mathrm{def}}{:=}\, 
 -\frac{2(1-16t) \mathscr{K}(16t)}{\mathscr{E}(16t)-(1-16t)\mathscr{K}(16t)}\left(\frac{2}{\pi}\, \mathscr{K}(16t)-1\right)
	\ln\left(1-\frac{\pi}{2\mathscr{K}(16t)}\right) + \ln(t) + 4 \ln(2)\, .
$$
\end{prop}

\begin{figure}
\begin{center}
\begin{tikzpicture}
	\draw (0,0) node[] {\includegraphics[height=7cm]{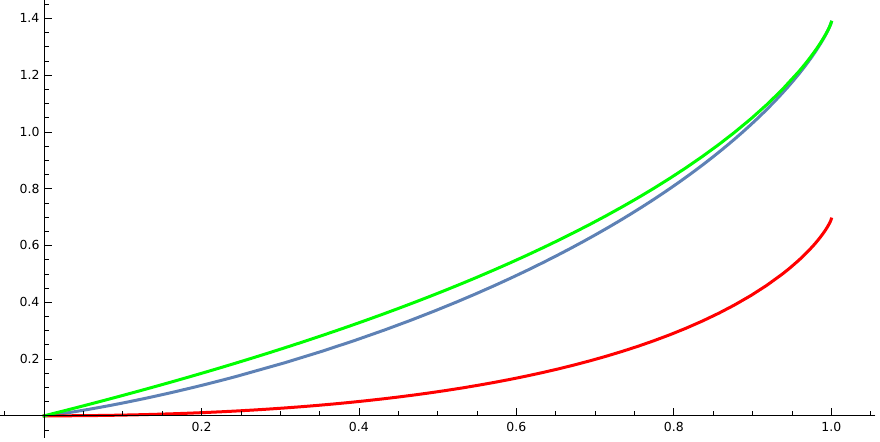}};
\end{tikzpicture}
	\caption{We plot all three large deviation rate functions together: blue is for the 2SRWB,
green is for the Dyck path, and red is for the SRWB. 
\label{fig:multipleRate}}
\end{center}
\end{figure}
In Figure \ref{fig:multipleRate}, we have plotted the rate function for the 2SRWB ensemble using parametric-plot,
along with the rate functions for the SRWB and Dyck ensembles.

\begin{remark}
As a consistency check, we note that 
$$
\left(\widetilde{X}(0),\widetilde{I}(0)\right)\, =\, (1,2\ln(2))\, ,
$$
which is consistent with equation (\ref{eq:right2SRWB}).
\end{remark}
\begin{remark}
We also have
$$
\left(\widetilde{X}\left(\frac{1}{16}\right),\widetilde{I}\left(\frac{1}{16}\right)\right)\, =\, (0,0)\, ,
$$
and
$$
	\frac{\widetilde{I}(\frac{1}{16}-\epsilon)}{\widetilde{X}(\frac{1}{16}-\epsilon)}\,
	=\, \frac{\pi}{\ln(1/\epsilon)}\
	\text{ and }\
	\widetilde{X}\left(\frac{1}{16}-\epsilon\right)\, \sim\, \epsilon\, \left(\ln(\epsilon)\right)^2\, .
$$
We note that for the 2SRWB problem, $\mathfrak{Z}_n$ is asymptotically distributed
as $\ln(n)\mathsf{T}/\pi$ where $\mathsf{T}$ is a Gamma-(2,1) random variable.
So, if we take  large but fixed $\tau$, then $\mathbf{P}(\mathfrak{Z}_n \geq \tau \ln(n))$
is asymptotically, for large $n$ equal to: $\mathbf{P}(\mathsf{T}\geq \pi \tau)
=\left(1+\pi \tau\right)\exp(-\pi \tau)$.
So $\frac{1}{n}\, \ln\left(\mathbf{P}(\mathfrak{Z}_n \geq \tau \ln(n)))\right) \sim -\pi \tau/n$.
So the large deviation rate function should agree with the asymptotics of the negative: $\pi \tau/n$ at its left endpoint.

Now take $\tau \ln(n) = t n$ so that $t = \tau \ln(n)/n$.
Note that this implies $\ln(t) \sim -\ln(n)$, as $n \to \infty$.
If the leading order asymptotics of $\widetilde{X}\left(\frac{1}{16}-\epsilon\right)$ -- namely $\epsilon (\ln(\epsilon))^2$ -- is set to be equal to $t$, then 
we have $\epsilon \sim t/(\ln(t))^2$.
That means $\ln(1/\epsilon) \sim \ln(1/t) \sim \ln(n)$ to leading order.
But notice that $t/(\ln(t))^2$ is in turn $\sim \tau/(n\ln(n))$, as $n \to \infty$.
So the large deviation rate gives 
$\widetilde{I}(\frac{1}{16}-\epsilon) \sim \pi \widetilde{X}\left(\frac{1}{16}-\epsilon\right)/\ln(1/\epsilon)$
is asymptotically $\pi \epsilon (\ln(\epsilon))^2/(\ln(1/\epsilon))
= \pi \epsilon \ln(1/\epsilon) \sim \pi\tau/n$.
This is consistent.
\end{remark}
The generating function for the combinatorial sequence $\boldsymbol{a}$ where $a_n = B_n^2$  for $|t|<1/16$ is
$$
	g(t)\, =\, \sum_{n=0}^{\infty} \binom{2n}{n}^2 t^n\, =\, \frac{2}{\pi}\, \mathscr{K}(16t)\, .
$$
Therefore, by Chebyshev's inequality, we may deduce that 
\begin{equation}
\label{ineq:Cheb2SRWB}
	\limsup_{\substack{n \to \infty\\k/n \to \kappa}} \frac{1}{n}\, \ln\left(\mathcal{M}_n^{(k)}(\boldsymbol{a})\right)\, 
\leq\, \min_{t \in [0,1/16]} \Big(\kappa \ln(g(t)) - \ln(t)\Big)\, .
\end{equation}
Note that the evaluation of the limit of $n^{-1} \ln(\mathcal{M}_n^{(k)})$ is the main step from the previous section: Lemma \ref{lem:LDP1}.
(Everything else is basically calculus formulas.)
One could try to use the general strategy for large deviations to prove that the right-hand-side of (\ref{ineq:Cheb2SRWB}) is also the liminf.
But instead we view a major tool to be an application of the circle method of Hardy and Ramanujan that we demonstrate
with the following result.
\begin{theorem}
\label{thm:circle}
For each $\kappa>0$, there is a unique arg-min of of the function $t\mapsto\kappa \ln(g(t)) - \ln(t)$.
If we call the arg-min  $T_*(\kappa)$, then 
the desired asymptotics are expressible in terms of this:
$$
	\mathcal{M}_n^{(k)}(\boldsymbol{a})\, \sim\, e^{k \ln(g(T_*(k/n)))-n \ln(T_*(k/n))}\, \frac{1}{\sqrt{2\pi k v(\kappa)}}\, ,
$$
as $n \to \infty$ assuming $k/n \to \kappa \in (0,\infty)$, where
$$
	v(\kappa)\, =\, t \cdot \frac{d}{dt}\, \ln(g(t)) + t^2\, \cdot \frac{d^2}{dt^2}\, \ln(g(t)) \Bigg|_{t=T_*(\kappa)}\, .
$$
\end{theorem}
This will be proved in Section \ref{sec:circle}.

Here the factor $1/\sqrt{2\pi k v(\kappa)}$ comes from the evaluation of an oscillatory integral over the circle,
using Hardy and Ramanujan's circle method.
In particular, all this does imply that 
\begin{equation}
\label{eq:mom2SRWB}
	\lim_{\substack{n \to \infty\\k/n \to \kappa}} \frac{1}{n}\, \ln\left(\mathcal{M}_n^{(k)}(\boldsymbol{a})\right)\, 
=\, \min_{t \in [0,1/16]} \Big(\kappa \ln(g(t)) - \ln(t)\Big)\, .
\end{equation}
That is the main input we need to prove Proposition \ref{prop:2SRWBrate}.

\section{Proofs of the ``first results''}
\label{sec:proofs}

Let us begin with the proof of the lemma.

Note that for each $\vec{x} \in \mathcal{X}_n$ and any $k \in \mathfrak{Z}_n(\vec{x})$, our axioms guarantee
that $(x_1,\dots,x_k)$ and $(x_{k+1},\dots,x_n)$ are elements of $\mathcal{X}_k$ and $\mathcal{X}_{n-k}$,
respectively (and also for $k=0$).
The axioms are such that the converse is also true.
Therefore, we may break up the sum as 
$$
\mathcal{M}_n^{(k)}(\boldsymbol{a})\, =\, \sum_{\vec{x} \in \mathcal{X}_n}
\sum_{j_1=0}^{\infty} \cdots \sum_{j_k=0}^{\infty} \mathbf{1}_{\{n\}}(j_1+\dots+j_k)
\prod_{r=1}^{k} \mathbf{1}_{\{\o\}}(x_{j_1+\dots+j_r})\, .
$$
But, writing $\ell_r=j_1+\dots+j_r$ for $r=1,\dots,n$, this sum is clearly the same as 
$$
\mathcal{M}_n^{(k)}(\boldsymbol{a})\, =\, \sum_{\vec{x} \in \mathcal{X}_n}
\left|\left\{(\ell_1,\dots,\ell_k) \in [n]^k\, :\, 0\leq \ell_1\leq \dots\leq \ell_k\leq n\ \text{ and }\ 
\ell_1,\dots,\ell_k \in \mathfrak{Z}_k(\vec{x})\right\}\right|\, .
$$
The summand is exactly the quantity being counted in Feller's ``stars and bars'' problem.
See page 38 of his \cite{Feller}.

We will not reproduce the proof of Proposition \ref{prop:Catalan}.
Please see the excellent article by Larcombe and French \cite{LarcombeFrench}.

\begin{proofof}{\bf Proof of Corollary \ref{cor:Dyck}:}
The uniform distribution on Dyck paths assigns probability $1/C_n$ to each valid Dyck path.
Since that is $\mathcal{M}^{(1)}_n(\boldsymbol{C})$, we have
$$
	\mathbf{E}^{\mathrm{Dyck}}\left[\binom{|\mathfrak{Z}_n|+k-1}{k}\right]\, =\, \frac{\mathcal{M}^{(k+1)}_n(\boldsymbol{C})}
{\mathcal{M}^{(1)}_n(\boldsymbol{C})}\, .
$$
Moreover, by Proposition \ref{prop:Catalan}, the right-hand-side is equal to 
$$
\frac{k+1}{2n+k+1} \cdot \frac{2n+1}{1}\, \cdot \frac{\binom{2n+1}{n}}{\binom{2n+k+1}{n}}\, 
\sim\, (k+1) 2^{-k}\, ,\ \text{ as $n \to \infty$,}
$$
where we use Stirling's formula to calculate the asymptotic expression.
Then note that by Newton's binomial formula (and the fact that $\binom{-x}{k} = (-1)^k \binom{x+k-1}{k}$)
we have
$$
\lim_{n \to \infty}
\mathbf{E}^{\mathrm{Dyck}}\left[\left(\frac{1}{1-t}\right)^{|\mathfrak{Z}_n|}\right]\,
=\, \frac{1}{(1-2t)^2}\, .
$$
Then changing variables to $z = 1/(1-t)$ and $t = (z-1)/z$, we obtain the probability generating function
$(z/2)^2 (1-(z/2))^{-2}$, which yields the stated probabilitie by extracting the coefficients
(by Newton's binomial formula again).
\end{proofof}

The proof of Proposition \ref{prop:central} is definitely simpler than the proof of Proposition \ref{prop:Catalan}
which we skipped. Let us note (using $(-1)^n\binom{-1/2}{n} = \binom{n-(1/2)}{n} = 2^{-2n} \binom{2n}{n}$) that since
$$
	\sum_{n=0}^{\infty} \binom{2n}{n} t^n\, =\,  \frac{1}{(1-4t)^{1/2}}\, ,
$$
we have $(\sum_{n=0}^{\infty} \binom{2n}{n} t^n)^k$ is equal $(1-4t)^{-k/2}$ which may be easily
expanded by Newton's binomial formula.
Therefore, since the generating function of the convolution is the power of the generating function,
by Cauchy's product formula, we obtain that formula.

\begin{proofof}{\bf Proof of Corollary \ref{cor:SRWB}:}
Note that
$$
	\mathbf{E}^{\mathrm{SRWB}}\left[\binom{|\mathfrak{Z}_n|+k-1}{k}\right]\, =\, \frac{\mathcal{M}^{(k+1)}_n(\boldsymbol{B})}
{\mathcal{M}^{(1)}_n(\boldsymbol{B})}\, .
$$
But by Proposition \ref{prop:central} this is 
$$
 \frac{\mathcal{M}^{(k+1)}_n(\boldsymbol{B})}
{\mathcal{M}^{(1)}_n(\boldsymbol{B})}\, =\,
\frac{(n-1+((k+1)/2))_n}{(n-1+(1/2))_ n}\, .
$$
Note that the well-known formula for the Gamma function (aside from its basic definition) is
$$
	\Gamma(x)\, =\, \lim_{n \to \infty} \frac{n!}{(n+s-1)_n}\, n^{s-1}\, .
$$
So we obtain
$$
\lim_{n \to \infty}  n^{-k/2} \frac{\mathcal{M}^{(k+1)}_n(\boldsymbol{B})}
{\mathcal{M}^{(1)}_n(\boldsymbol{B})}\, =\,
\lim_{n\to \infty} \frac{n!}{(n-1+(1/2))_ n}\, n^{(1/2)-1} 
\left( \frac{n!}{(n-1+((k+1)/2))_n}\, n^{((k+1)/2)-1}\right)^{-1}\, ,
$$
which equals $\Gamma\left(\frac{1}{2}\right) / \Gamma\left(\frac{k+1}{2}\right)$.
Note that by Legendre's duplication formula this is $2^k \Gamma\left(\frac{k}{2}+1\right)/\Gamma(k+1)$.
But 
$$
\lim_{n \to \infty} \mathbf{E}^{\mathrm{SRWB}}\left[n^{-k/2} \binom{|\mathfrak{Z}_n|+k-1}{k}\right]\, 
=\, 
\lim_{n \to \infty} \mathbf{E}^{\mathrm{SRWB}}\left[\frac{1}{k!} \left(n^{-1/2} |\mathfrak{Z}_n|\right)^k\right]\, 
$$
as may be seen by induction. Since $\Gamma(k+1) = k!$, this yields the desired result.
\end{proofof}
For the proofs of Theorems \ref{thm:2SRWB} and \ref{thm:3SRWB}, we will just prove the asymptotic
formulas for the convolutions. The probabilistic consequences then follow as in the proofs of 
Corollary \ref{cor:SRWB} and \ref{cor:Dyck}, respectively.
\begin{proofof}{\bf Proof of Theorem \ref{thm:2SRWB}:}
We have, by Stirling's formula,
$$
	\binom{2n}{n}^2\, \sim\, 2^{4n\ln(2)}\, \frac{1}{\pi n}\, ,\ \text{ as $n \to \infty$.}
$$
Now we note that
$$
	\sum_{k=1}^{n-1} \frac{1}{k(n-k)}\, =\, \frac{1}{n}\, \sum_{k=1}^{n-1} \left(\frac{1}{k} + \frac{1}{n-k}\right)\,
\sim\, \frac{2\ln(n)}{n}\, ,\ 
\text{ as $n \to \infty$,}
$$
and
for $r \in \N$,
\begin{equation*}
\begin{split}
	\sum_{k=1}^{n-1} \frac{\left(\ln(k)\right)^r}{k(n-k)}\,
	&=\, \frac{1}{n}\, 
\sum_{k=1}^{n-1} \frac{\left(\ln(n) + \ln(k/n)\right)^r}{(k/n)(1 - \frac{k}{n})}  \cdot \frac{1}{n}\\
	&\sim\, \frac{2\left(\ln(n)\right)^{r+1}}{n} 
+ \frac{1}{n}\, \sum_{s=1}^{r} \binom{r}{s} \left(\ln(n)\right)^{r-s} \int_{1/n}^{1-(1/n)} 
\frac{\left(\ln(t)\right)^s}{t(1-t)}\, dt\\
	&\sim\, \frac{2\left(\ln(n)\right)^{r+1}}{n}
+\frac{1}{n}\, \sum_{s=1}^{r} \binom{r}{s} (-1)^s \frac{\left(\ln(n)\right)^{r+1}}{s+1}\\
	&=\, \frac{(r+2) \left(\ln(n)\right)^{r+1}}{(r+1)n}\, .
\end{split}
\end{equation*}
Using this, one can prove the formula by induction:
$$
\sum_{1<j_1<j_2<\dots<j_{r+1}<n} \frac{1}{j_1} \cdot \frac{1}{j_2-j_1} \cdots \frac{1}{j_{r+1}-j_r} \cdot \frac{1}{n-j_{r+1}}\, \sim\, \frac{(r+2) \left(\ln(n)\right)^{r+1}}{n}\, .
$$
This directly implies the result.
\end{proofof}
The difference between Theorem \ref{thm:2SRWB} and \ref{thm:3SRWB} is that in the latter the asymptotics
of the moment sequences does not need a recaling factor depending on $n$. 

\begin{proofof}{\bf Proof of Theorem \ref{thm:3SRWB}:}
We now use the formula that for each fixed $k \in \{0,1,\dots\}$, we have
$$
\frac{\binom{2n-2k}{n-k}^3}{\binom{2n}{n}^3}\, \sim\, 2^{-6k\ln(2)}\, ,\ \text{ as $n \to \infty$.}
$$
So then we can see that
$$
\frac{1}{\binom{2n}{n}^3}\, \sum_{1<j_1<j_2<\dots<j_{r}<n} 
\binom{2j_1}{j_1}^3 \binom{2(j_2-j_1)}{j_2-j_1}^3 \cdots
\binom{2(j_r-j_{r-1})}{j_{r}-j_{r-1}}^2 \cdot \binom{2(n-j_r)}{n-j_{r}}^3\, 
$$
is asymptotic to 
$$
	\sum_{s=0}^{r} \sum_{\ell_1=0}^\infty \cdots \sum_{\ell_s=0}^{\infty} \sum_{\ell_{s+1}=0}^\infty
\cdots \sum_{\ell_{r}=0}^{\infty} \prod_{k=1}^{r-1} \left(2^{-6 \ell_k} \binom{2\ell_k}{\ell_k}^3\right)\, ,
$$
where for a given choice of $s$, we let $\ell_k = j_k-j_{k-1}$ (setting $j_0=0$) for $k\in \{1,\dots,s\}$
and $\ell_{k} = j_{k+1}-j_{k}$ (setting $j_{r+1}=n$) for $k \in \{s+1,\dots,r\}$.
Note that then $\ell_1+\dots+\ell_r=n-(k_{s+1}-k_s)$, so that $k_{s+1}-k_s=n-(\ell_1+\dots+\ell_r)$
which is ``close to $n$'' as long as $\ell_1,\dots,\ell_r$ are ``order-1''.
So we can use the asymptotic formula above for $\binom{2(k_{s+1}-k_s)}{k_{s+1}-k_s}/\binom{2n}{n}$:
it is asymptotic to $2^{-6(\ell_1+\dots+\ell_r)}$.
(Note that at least 1 of the summands must be large.
But we may prove by induction that the contributions sum to a negligible
amount for those choices where more than 1 summand is large.)
Since $\sum_{n=0}^{\infty} \binom{2n}{n}^3 \frac{1}{64^n} = \pi/\left(\Gamma\left(\frac{3}{4}\right)\right)^3$,
we have the result.
\end{proofof}

\section{Proof of Theorem \ref{thm:rate} using calculus formulas}
\label{sec:rates}

Let us discuss the LDP with respect to the SRWB ensemble, since the Dyck path
ensemble works the same.
For a fixed $n$, decompose $x$ into subintervals $\mathcal{J}_n^{(i)} = [i/n,(i+1)/n]$ for $i=0,\dots,n-1$.
Note that $\mathbf{P}^{\mathrm{SRWB}}_n(\mathfrak{Z}_n \in \mathcal{J}_n^{(i)})$
is bounded above by $\mathbf{P}^{\mathrm{SRWB}}_n(\mathfrak{Z}_n \geq i/n)$
which is bounded above by $\exp\left(-nI^{\mathrm{SRWB}}(i/n)\right)$ by Chebyshev's inequality.
Fix a $y >0$. Then by uniform continuity, we know that $U(x,y)$ varies by at most $\delta_n$
on each $\mathcal{J}_n^{(i)}$, for $i \in \{0,\dots,n-1\}$, and $\lim_{n \to \infty} \delta_n  =0$.
Hence,
$$
	n^{-1}\, \ln\left(\mathbf{E}^{\mathrm{SRWB}}_n\left[e^{n U(\mathfrak{Z}_n/n,y)} \mathbf{1}_{\mathcal{J}_n^{(i)}}(\mathfrak{Z}_n/n)\right]\right)\, \leq\, W\left(\frac{i}{n},y\right)+\delta_n\, .
$$
Choose $\epsilon>0$. Then for any $i \in \{0,\dots,n-1\}$ such that $|(i/n)-X_*^{\mathrm{SRWB}}(y)|>\epsilon$, we have $W(i/n,y)<\mathcal{U}^{\mathrm{SRWB}}(y)-C(y) \epsilon^{1/2}$ by strict concavity
of $W(\cdot,y)$, where $C(y)$ is bounded blow by the second derivative of $W(\cdot,y)$
at $X^{\mathrm{SWRB}}_{*}(y)$.
So we have
$$
	n^{-1}\, \ln\left(\mathbf{E}^{\mathrm{SRWB}}_n\left[e^{n U(\mathfrak{Z}_n/n,y)} 
\mathbf{1}_{(X^{\mathrm{SRWB}}_*(y)-\epsilon,X^{\mathrm{SRWB}}_*(y)+\epsilon)^c}(\mathfrak{Z}_n/n)\right]\right)\, \leq\, \mathcal{U}^{\mathrm{SRWB}}(y)+\delta_n - C(y)\epsilon^{1/2} + n^{-1}\ln(n)\, .
$$
Hence,
$$
	\limsup_{n \to \infty} n^{-1}\, \ln\left(\mathbf{E}^{\mathrm{SRWB}}_n\left[e^{n U(\mathfrak{Z}_n/n,y)} 
\mathbf{1}_{(X^{\mathrm{SRWB}}_*(y)-\epsilon,X^{\mathrm{SRWB}}_*(y)+\epsilon)^c}(\mathfrak{Z}_n/n)\right]\right)\, \leq\, \mathcal{U}^{\mathrm{SRWB}}(y)- C(y)\epsilon^{1/2}\, .
$$
But since 
$$
	\lim_{n \to \infty} n^{-1}\, \ln\left(\mathbf{E}^{\mathrm{SRWB}}_n\left[e^{n U(\mathfrak{Z}_n/n,y)} 
\right]\right)\, \leq\, \mathcal{U}^{\mathrm{SRWB}}(y)\, ,
$$
that means
$$
	\liminf_{n \to \infty} n^{-1}\, \ln\left(\mathbf{E}^{\mathrm{SRWB}}_n\left[e^{n U(\mathfrak{Z}_n/n,y)} 
\mathbf{1}_{[X^{\mathrm{SRWB}}_*(y)-\epsilon,X^{\mathrm{SRWB}}_*(y)+\epsilon]}(\mathfrak{Z}_n/n)\right]\right)\, \geq\, \mathcal{U}^{\mathrm{SRWB}}(y)\, .
$$
By uniform continuity of $U(\cdot,y)$ again, we have that it varies by at most $\widetilde{\delta}(\epsilon)$
on $[X^{\mathrm{SRWB}}_*(y)-\epsilon,X^{\mathrm{SRWB}}_*(y)+\epsilon]$ for a function
such that $\lim_{\epsilon \to 0^+} \widetilde{\delta}(\epsilon)=0$.
From this we can see that
$$
	\liminf_{n \to \infty} n^{-1}\, \ln\left(\mathbf{E}^{\mathrm{SRWB}}_n\left[ 
\mathbf{1}_{[X^{\mathrm{SRWB}}_*(y)-\epsilon,X^{\mathrm{SRWB}}_*(y)+\epsilon]}(\mathfrak{Z}_n/n)\right]\right)\, \geq\, \mathcal{U}^{\mathrm{SRWB}}(y)-U(X^{\mathrm{SRWB}}_*(y),y)-\widetilde{\delta}(\epsilon)\, .
$$
And this means
$$
	\liminf_{n \to \infty} n^{-1}\, \ln\left(\mathbf{P}^{\mathrm{SRWB}}_n\left(
n^{-1} \mathfrak{Z}_n \in [X^{\mathrm{SRWB}}_*(y)-\epsilon,X^{\mathrm{SRWB}}_*(y)+\epsilon]\right)\right)\, \geq\, -I^{\mathrm{SRWB}}(X^{\mathrm{SRWB}}_*(y))-\widetilde{\delta}(\epsilon)\, .
$$
In particular, this means
$$
	\liminf_{n \to \infty} n^{-1}\, \ln\left(\mathbf{P}^{\mathrm{SRWB}}_n\left(
n^{-1} \mathfrak{Z}_n \geq X^{\mathrm{SRWB}}_*(y)-\epsilon\right)\right)\, \geq\, -I^{\mathrm{SRWB}}(X^{\mathrm{SRWB}}_*(y))-\widetilde{\delta}(\epsilon)\, .
$$
Taking the limit $\epsilon\to0^+$ gives the bound:
$$
	\lim_{\epsilon \to 0^+} \liminf_{n \to \infty} n^{-1}\, \ln\left(\mathbf{P}^{\mathrm{SRWB}}_n\left(
n^{-1} \mathfrak{Z}_n \geq X^{\mathrm{SRWB}}_*(y)-\epsilon\right)\right)\, \geq\, -I^{\mathrm{SRWB}}(X^{\mathrm{SRWB}}_*(y))\, .
$$
But the mapping $y \to X^{\mathrm{SRWB}}_*(y)$ is invertible: its inverse is $x \mapsto Y^{\mathrm{SRWB}}_*(x)$.
So, in particular, it is onto. Hence, for all $x \in (0,1)$ we have
$$
	\lim_{\epsilon \to 0^+} \liminf_{n \to \infty} n^{-1}\, \ln\left(\mathbf{P}^{\mathrm{SRWB}}_n\left(
n^{-1} \mathfrak{Z}_n \geq x-\epsilon\right)\right)\, \geq\, -I^{\mathrm{SRWB}}(x)\, .
$$
For any $x \in (0,1)$, we may now take $x'$ values slightly greater than $x$, use monotonicity of probability, and use continuity of $I^{\mathrm{SRWB}}$
to conclude
$$
	\liminf_{n \to \infty} n^{-1}\, \ln\left(\mathbf{P}^{\mathrm{SRWB}}_n\left(
n^{-1} \mathfrak{Z}_n \geq x\right)\right)\, \geq\, -I^{\mathrm{SRWB}}(x)\, .
$$

\section{Proof of calculus formula for Theorem \ref{thm:rate}}
\label{sec:LEMrates}

\begin{proofof}{\bf Proof of Lemma \ref{lem:LDP1}:}
We note that
$$
	n^{-1}\, \ln \left(\binom{\mathfrak{Z}_n+k-1}{k}\right)\, \sim\, U\left(\frac{\mathfrak{Z}_n}{n}\, ,\ 
\frac{k}{n}\right)$$
Then the formulas follow from the precise asymptotics of the well-known formulas in Proposition \ref{prop:Catalan} and  Proposition \ref{prop:central}.
\end{proofof}

\begin{proofof}{\bf Proof of Lemma \ref{lem:CalFacts1}:}
We first note that
\begin{equation}
\label{eq:VFirstDeriv}
	\frac{\partial}{\partial y}\, V^{\mathrm{SRWB}}(x,y)\, =\, \ln\left(\frac{\sqrt{y(y+2)}}{x+y}\right)\
\text{ and }\
	\frac{\partial}{\partial y}\, V^{\mathrm{Dyck}}(x,y)\, =\, \ln\left(\frac{y(2+y)}{(1+y)(x+y)}\right)\, .
\end{equation}
It is easy to see that the derivatives are $0$ if $y$ is equal to $Y^{\mathrm{SRWB}}_*(x)$
and $Y^{\mathrm{Dyck}}_*(x)$, respectively.
(The arguments of the logarithm are equal to $1$.)
Moreover
$$
	\frac{y(2+y)}{(1+y)(x+y)}\Bigg|_{y=Y_*^{\mathrm{Dyck}}(x)+\alpha}\,
	=\, \frac{a + b \alpha}{a+c \alpha}\, ,\ \text{ for }\
	a\, =\, \frac{x(2-x)}{(1-x)^2} + \alpha^2\, ,\ 
	b\, =\, \frac{2}{(1-x)^2}\ \text{ and }\ c\, =\, \frac{1+2x-x^2}{(1-x)^2}\, .
$$
Since $a>0$ and $b>c$, and $a+b\alpha>0$, $a+c\alpha>0$ (for $x \in (0,1)$ and $y>0$) we see that for $\alpha>0$ the $y$-partial derivative of $V^{\mathrm{Dyck}}$ is positive
and for $\alpha<0$ it is negative. Hence $-V^{\mathrm{Dyck}}(x,\cdot)$ is unimodular.
In other words, $V^{\mathrm{Dyck}}(x,\cdot)$ is quasi-convex.
It has a unique minimizer.
A similar calculation shows that $V^{\mathrm{SRWB}}(x,\cdot)$ has the same property
because
$$
	\frac{\sqrt{y(y+2)}}{x+y}\, =\, \frac{\sqrt{(x+y)^2 + 2(1-x)\alpha}}{x+y}\, ,\ \text{ for }
	\alpha\, =\, y - Y^{\mathrm{SRWB}}_*(x)\, .
$$
The formulas in (\ref{eq:IRateSub}) are direct substitutions, which may be verified
with a symbolic algebra program such as Wolfram Mathematica.
\end{proofof}

\begin{proofof}{\bf Proof of Lemma \ref{lem:CalFacts2}:}
Direct calculation shows
\begin{equation}
\label{eq:WFirstDeriv}
	\frac{\partial}{\partial x}\, W^{\mathrm{SRWB}}(x,y)\, =\, \ln\left(\frac{(1-x)(x+y)}{x(2-x)}\right)\
\text{ and }\
	\frac{\partial}{\partial x}\, W^{\mathrm{Dyck}}(x,y)\, =\, \ln\left(\frac{2(1-x)(x+y)}{x(2-x)}\right)\, .
\end{equation}
Hence,
$$
	\frac{\partial^2}{\partial x^2}\, W^{\mathrm{SRWB}}(x,y)\, =\,
	\frac{\partial^2}{\partial x^2}\, W^{\mathrm{Dyck}}(x,y)\, =\, 
	-\frac{(1+y)x^2-2xy+2y}{x(1-x)(2-x)(x+y)}\, .
$$
This may be rewritten as 
$$
	-\frac{((1+y)x-y)^2+y(y+2)}{x(1-x)(2-x)(1+y)(x+y)}\, .
$$
Hence the functions are strictly concave.
It is easy to see that the conditions to get 0 for the derivatives in (\ref{eq:WFirstDeriv})
are as stated.
Then equation (\ref{eq:Wsub}) follows from direct substitution which may be simplified
in a symbolic algebra program such as Wolfram Mathematica.
\end{proofof}

\section{A simplified circle method}

\label{sec:circle}

We now turn to the proof of Theorem \ref{thm:circle}.
We believe that the result holds in some generality. The reason is that $k$ is large if $k \sim \kappa n$ for $\kappa \in (0,\infty)$ fixed, as $n \to \infty$.
By Cauchy's integral formula, we know that
$$
	\mathcal{M}_n^{(k)}(\boldsymbol{a})\, 
	=\, \frac{1}{t^n}\, \int_{-\pi}^{\pi} 
	\left(g\left(t e^{i\theta}\right)\right)^k\,
	e^{-in\theta}\, \frac{d\theta}{2\pi}\, ,
$$
for every $t>0$ such that $t$ does not exceed the radius of convergence of $g$.
We believe that the usual Hayman method should be applicable. Considering the major arc, one does obtain the stated asymptotics.
We note that a similar formula that holds in great generality is the Bahadur-Rao theorem \cite{BahadurRao}.
We cite Dembo and Zeitouni for a pedagogical proof of that formula \cite{DemboZeitouni}, Section 3.7.

The only difficulty is in proving the usual bounds in the minor arc. 
However, since we do not yet have a general proof for the minor arc, we will instead consider the theorem just for the special case of the 2SRWB model, which is our primary interest, here.
And we will use a special trick, which is just the application of the Laplace transform, since we happen to know how to obtain the complete elliptic integral 
of the first kind using a Laplace transform.

The Hardy and Ramanujan circle method was first established in \cite{HardyRamanujan}.
There were many simplifications of the method over the years.
We mention work of Rademacher \cite{Rademacher}, Newman \cite{Newman} for Hardy and Ramanujan's specific example of the partition numbers.
For general applications, we cite the textbook of Flajolet and Sedgewick \cite{FlajoletSedgewick}
as well as more recent work by Cant\'on, et.~al, \cite{CantonFFM}. (We thank Malek Abdesselam for bringing this reference to our attention.)
But the simplification we are going to most closely follow is one started by Fristedt \cite{Fristedt} and Romik \cite{Romik}.
We first note that the generating function may be written as 
$$
g(t)\, =\, \frac{2}{\pi}\, \mathbf{E}\left[\frac{1}{\sqrt{1-\mathsf{U}^2}}\, \cdot 
\frac{1}{\sqrt{1-16t\mathsf{U}^2}}\right]\, ,
$$
where $\mathsf{U}$ is a continuous uniform random variable in $(0,1)$, because this is the same
as
$$
\frac{2}{\pi}\, \mathscr{K}(16t)\ \text{ for }\
\mathscr{K}(z)\, =\, \int_0^1 \frac{1}{\sqrt{(1-x^2)(1-zx^2)}}\, dx\, .
$$
Then, using Gamma-$(1/2,1)$ random variables, we may rewrite this as 
$$
g(t)\, =\, \frac{2}{\pi}\, \mathbf{E}\left[e^{\mathsf{U}^2 \mathsf{X}} e^{16t\mathsf{U}^2 \mathsf{Y}}\right]\, ,
$$
where (independently of $\mathsf{U}$) we have that $\mathsf{X}$ and $\mathsf{Y}$ are independent Gamma-$(1/2,1)$ random variables with 
probability density function
$$
f(x)\, =\, \frac{x^{-1/2} e^{-x}}{\Gamma(1/2)}\, ,\ \text{ for $x>0$.}
$$
Using random variables to rewrite the generating function is a convenience,
but one whose inspiration we believe may be traced to Fristedt and Romik.
We feel the key step is the introduction of the Laplace transform,
which this does in passing.
The Laplace transform is not part of Hayman's method.
But it does show up in the multiplicative version of the Mellin transform frequently in analytic number theory.
It is really a well-established method. We note that for extracting asymptotics of algebraically singular generating functions (at their radius of convergence)
one does need some extra tool such as the Laplace method.
For example if one wished to extract asymptotics of the $k$-fold convolution of any of the sequences considered here, by the saddle-point method,
for finite $k$, then one would need the Laplace transform. (Of course, since the formulas for the convolutions may be determined more simply by methods
other than the saddle-point method, that is not necessary.)
We cite a pedagogical survey of one of the authors \cite{Starr} which will remain in preprint form (since it merely re-establishes well-known results in a more pedantic way).
As we indicated before, since $k$ is large in Theorem \ref{thm:circle}, the Laplace transform trick is not strictly necessary.
But we do find it convenient in the present framework.

We note that if we define the tilted measure $\mu_t$ on $(0,\infty)$ by
$$
\mu_t(A)\, =\, \frac{\mathbf{E}\left[\mathbf{1}_A(16\mathsf{U}^2\mathsf{Y})
e^{\mathsf{U}^2 \mathsf{X}} e^{16t\mathsf{U}^2 \mathsf{Y}}\right]}
{\mathbf{E}\left[e^{\mathsf{U}^2 \mathsf{X}} e^{16t\mathsf{U}^2 \mathsf{Y}}\right]}\, ,
$$
then
$$
\frac{d}{dt}\, \ln(g(t))\, =\, \int_0^{\infty} x\, d\mu_t(x)\ \text{ and }\ 
\frac{d^2}{dt^2}\, \ln(g(t))\, =\, \int_0^{\infty} x^2\, d\mu_t(x) - \left(\int_0^{\infty} x\, d\mu_t(x)\right)^2\, .
$$
These are $\operatorname{avg}(\mu_t)$ and $\operatorname{var}(\mu_t)$, respectively.
So clearly $\ln(g(t))$ is strictly convex.
Thus we see that
$$
\frac{d}{dt}\, \left(\kappa \ln(g(t))-\ln(t)\right)\,
	=\, \kappa \operatorname{avg}(\mu_t) - \frac{1}{t}\, ,
$$
and
$$
\frac{d^2}{dt^2}\, \left(\kappa \ln(g(t))-\ln(t)\right)\,
	=\, \kappa \operatorname{var}(\mu_t) + \frac{1}{t^2}\, .
$$
Hence, $\kappa \ln(g(t))-\ln(t)$ is strictly convex on $(0,1/16)$.
Since $g(0)=1$ and $\lim_{t \to 1^-} g(t)=\infty$, we know that there is a local minimizer of 
$\kappa \ln(g(t))-\ln(t)$.
By strict convexity it is a unique, global minimizer.
We call it $T_*(\kappa)$. 
Also, the critical point equation, which is satisfied by the arg-min $T_*(k)$ is that
\begin{equation}
\label{eq:Cp}
\Bigg(\frac{d}{dt}\, \left(\kappa \ln(g(t))-\ln(t)\right)\bigg|_{t=T_*(\kappa)}\, =\, 0\Bigg)\qquad \Leftrightarrow\qquad
\Bigg(\kappa T_*(\kappa) \operatorname{avg}(\mu_{T_*(\kappa)})\, =\, 1\Bigg)\, .
\end{equation}


So, by the Cauchy integral formula, taking independent copies of the random variables
we have
$$
	\mathcal{M}_n^{(k)}(\boldsymbol{a})\, =\, \frac{(2/\pi)^k}{t^n}\,
	\int_{-\pi}^{\pi} \mathbf{E}\left[\exp\left(\sum_{j=1}^{k} \mathsf{U}_j^2 \mathsf{X}_j\right)
\exp\left(16 te^{i\theta}\sum_{j=1}^{k} \mathsf{U}_j^2 \mathsf{Y}_j\right)\right]
e^{-in\theta}\, \frac{d\theta}{2\pi}\, .
$$
We prefer to factorize this as 
\begin{equation*}
\begin{split}
	\mathcal{M}_n^{(k)}(\boldsymbol{a})\, &=\, \frac{(2/\pi)^k}{t^n}\,
	 \mathbf{E}\left[\exp\left(\sum_{j=1}^{k} \mathsf{U}_j^2 \mathsf{X}_j\right)
\exp\left(16 t\sum_{j=1}^{k} \mathsf{U}_j^2 \mathsf{Y}_j\right)\right]\\
&\hspace{1cm}
\int_{-\pi}^{\pi} \mathbf{E}\left[\exp\left(t(e^{i\theta}-1)\sum_{j=1}^{k} \mathsf{Z}_j(t)\right)\right]\,
e^{-in\theta}\, \frac{d\theta}{2\pi}\, ,
\end{split}
\end{equation*}
where $\mathsf{Z}_1(t),\dots,\mathsf{Z}_k(t)$ are IID, $\mu_t$-distributed random variables.
So this may be rewritten as 
$$
	\mathcal{M}_n^{(k)}(\boldsymbol{a})\, =\, \frac{\left(g(t)\right)^k}{t^n}\,
\int_{-\pi}^{\pi} \mathbf{E}\left[\exp\left( t(e^{i\theta}-1)\sum_{j=1}^{k} \mathsf{Z}_j(t)\right)\right]\,
e^{-in\theta}\, \frac{d\theta}{2\pi}\, ,
$$
Let us denote the circle integral as 
$$
	\mathcal{C}_{n,k}(t)\, =\, \int_{-\pi}^{\pi} \mathbf{E}\left[\exp\left( t(e^{i\theta}-1)\sum_{j=1}^{k} \mathsf{Z}_j(t)\right)\right]\,
e^{-in\theta}\, \frac{d\theta}{2\pi}\, .
$$
We note that by (\ref{eq:Cp}), we have
$$
	\mathcal{C}_{n,k}(T_*(k/n))\, =\, \left.\int_{-\pi}^{\pi} \mathbf{E}\left[\exp\left( t(e^{i\theta}-1)\sum_{j=1}^{k} \mathsf{Z}_j(t)\right)\right]\,
e^{-ikt\theta \operatorname{avg}(\mu_t) }\, \frac{d\theta}{2\pi}\, \right|_{t=T_*(k/n)}\, .
$$
We may rewrite $e^{i\theta}-1$ as $i\sin(\theta) -2\sin^2(\theta/2)$. So
$$
	\mathcal{C}_{n,k}(T_*(k/n))\, =\, \int_{-\pi}^{\pi} \mathbf{E}\left[e^{it\sum_{j=1}^{k} \left(\sin(\theta)\mathsf{Z}_j(t)-\theta \operatorname{avg}(\mu_t)\right)}
e^{-2t\sin^2(\theta/2)\sum_{j=1}^{k} \mathsf{Z}_j(t)}\right]\, \frac{d\theta}{2\pi}\, \Bigg|_{t=T_*(k/n)}\, .
$$
For the first factor inside the expectation we use the central limit theorem.
For the second factor we use the mean-value and Cram\'er's large deviation principle for $\mu_t$.
Note that for fixed $t<1/16$, the radius of convergence of the moment generating function
for $\mu_t$ is $(1/16)-t$. So Cram\'er's theorem does apply to the IID sum $\mathsf{Z}_1(t)+\dots+\mathsf{Z}_k(t)$, with a rate function: call it $J_t$.

\begin{lemma}
With probability 1 we have
$$
\left|\exp\left( t(e^{i\theta}-1)\sum_{j=1}^{k} \mathsf{Z}_j(t)\right)\right|\, \leq\, 1\, .
$$
More precisely, for $\theta \in (-\pi,\pi)$ we have, using $\cos(\theta)-1\leq -2\theta^2/\pi^2$,
$$
\left|\exp\left( t(e^{i\theta}-1)\sum_{j=1}^{k} \mathsf{Z}_j(t)\right)\right|\, \leq\, 
\exp\left(-2 t(\theta/\pi)^2\sum_{j=1}^{k} \mathsf{Z}_j(t)\right)\, .\ \square
$$
\end{lemma}
Let $\mathcal{A}_k(t,\epsilon)$ be the event that $\frac{1}{k} \sum_{j=1}^{k} \mathsf{Z}_j(t)$ is greater than or equal to $\operatorname{avg}(\mu_t)-\epsilon$.
Let $p_k(t,\epsilon) = 1-\mathbf{P}(\mathcal{A}_k(t,\epsilon))$ so that 
$$
	\lim_{k \to \infty} k^{-1} \ln(p_k(t,\epsilon))\, =\, - J_t(-\epsilon)\, <\, 0\, ,
$$
by Cram\'er's theorem.
Then, putting this together, we may conclude that 
\begin{equation*}
\begin{split}
	\mathcal{C}_{n,k}(T_*(k/n))\, &=\, \int_{-\ln(k)/\sqrt{k}}^{\ln(k)/\sqrt{k}} \mathbf{E}\left[e^{it\sum_{j=1}^{k} \left(\sin(\theta)\mathsf{Z}_j(t)-\theta \operatorname{avg}(\mu_t)\right)}
e^{-2t\sin^2(\theta/2)\sum_{j=1}^{k} \mathsf{Z}_j(t)}\right]\, \frac{d\theta}{2\pi}\, \Bigg|_{t=T_*(k/n)}\\
&\qquad
+ p_k(T_*(k/n),\epsilon) \cdot \mathcal{R}+ o\left(\frac{1}{\sqrt{k}}\right)\, ,
\end{split}
\end{equation*}
where $-1\leq \mathcal{R}\leq 1$.
Then we may conclude that 
$$
e^{-\delta(n,k)} \widetilde{\mathcal{C}}_{n,k}(-\epsilon) -p_k(T_*(k/n),\epsilon) + o\left(\frac{1}{\sqrt{k}}\right)\, 
	\leq\, \mathcal{C}_{n,k}(T_*(k/n))\, \leq\, e^{\delta(n,k)} \widetilde{\mathcal{C}}_{n,k}(\epsilon) + p_k(T_*(k/n),\epsilon) + o\left(\frac{1}{\sqrt{k}}\right)\, ,
$$
where
$$
\widetilde{\mathcal{C}}_{n,k}(\alpha)\, 
=\, \int_{-\ln(k)/\sqrt{k}}^{\ln(k)/\sqrt{k}} \mathbf{E}\left[e^{it\theta\sum_{j=1}^{k} \left(\mathsf{Z}_j(t)-\operatorname{avg}(\mu_t)\right)}\right]\,
e^{-tk\theta^2(\operatorname{avg}(\mu_t)+\alpha)/2}\, \frac{d\theta}{2\pi}\, \Bigg|_{t=T_*(k/n)}\, .
$$
The quantity $\delta(n,k)$ converges to $0$ as $n \to \infty$ with $k \sim \kappa n$, and is due to Taylor expanding $\sin(\theta)$ and $\sin^2(\theta/2)$.
But now we may use the CLT to conclude that
$$
\widetilde{\mathcal{C}}_{n,k}(\alpha)\, 
\sim\, \int_{-\ln(k)/\sqrt{k}}^{\ln(k)/\sqrt{k}} \mathbf{E}\left[e^{it\theta \sqrt{k\cdot \operatorname{var}(\mu_t)}\, \mathsf{N}}\right]\,
e^{-tk\theta^2(\operatorname{avg}(\mu_t)+\alpha)/2}\, \frac{d\theta}{2\pi}\, \Bigg|_{t=T_*(k/n)}\, ,
$$
where $\mathsf{N}$ is $\mathcal{N}(0,1)$-distributed.
So, then using the characteristic function of a standard-normal random variable, we have
$$
\widetilde{\mathcal{C}}_{n,k}(\alpha)\, 
\sim\, \int_{-\ln(k)/\sqrt{k}}^{\ln(k)/\sqrt{k}} e^{-t^2\theta^2k\cdot \operatorname{var}(\mu_t)/2}
e^{-tk\theta^2(\operatorname{avg}(\mu_t)+\alpha)/2}\, \frac{d\theta}{2\pi}\, \Bigg|_{t=T_*(k/n)}\, .
$$
Then the result follows by the usual Gaussian integral, for $\theta$, as in the usual Hayman method (as in \cite{FlajoletSedgewick}).
We may take $\epsilon \to 0^+$ after taking the $n \to \infty$ limit with $k \sim n \kappa$.
\section{A general derivation of the explicit rate function}

Here we include a section not directly germaine to the proofs of the theorems.
But its inclusion will assist in the proof of Proposition \ref{prop:2SRWBrate}.

We argue that if we have the explicit formulation for the asymptotics of $\mathcal{M}_{k}^{(n)}(\boldsymbol{a})$,
then we can determine the explicit rate function.

The main idea of the proof is similar to the Feynman-Hellman exercise:
because an optimizer satisfies a certain derivative equals 0,
by Fermat's lemma, we know that some other calculation works out more easily
than it otherwise would.
The Feynman-Hellman exercise applies this idea to first order perturbation
theory for eigenvalues of self-adjoint matrices, using the Rayleigh-Ritz
variational method. But the underlying idea
in optimization theory is not unique
to that setting.

By Laplace's method or Varadhan's lemma we have
\begin{equation*}
\max_{x \in [0,1]} \left(-I(x) + (x+\kappa) \ln(x+\kappa) - x \ln(x) - \kappa \ln(\kappa)
\right)\,
=\, \lim_{n \to \infty} \frac{1}{n}\, 
\left(\ln\left(\mathcal{M}_n^{(k)}(\boldsymbol{a})\right)- \ln(a_n)\right)\, ,
\end{equation*}
for $\boldsymbol{a}$ either equal to $\boldsymbol{B}$ or $\boldsymbol{C}$
and $I$ equal to $I^{\mathrm{SRWB}}$ or $I^{\mathrm{Dyck}}$, respectively.
Here we assume that $k$ depends on $n$ in such a way that $k/n \to \kappa \in (0,\infty)$ as $n \to \infty$.

Let us write $\kappa$ as $y$ and note that for the special value
$$
x = K(y)\, ,
$$
that optimizes the left-hand-side we should have by Fermat's lemma
$$
I'(K(y))\, =\, \frac{\partial}{\partial x}\left(
(x+y) \ln(x+y) - x \ln(x) - y \ln(y)
\right) \Bigg|_{x=K(y)}\, .
$$
Let us also define $Y = K^{-1}$. So then we have
\begin{equation*}
\begin{split}
&\hspace{-1cm}
\frac{d}{dx}\, \left(
(x+Y(x)) \ln(x+Y(x)) - x \ln(x) 
- Y(x) \ln(Y(x))
\right)\\
&\hspace{1cm}
=\, \frac{\partial}{\partial x}\left(
(x+y) \ln(x+y) - x \ln(x) - y \ln(y)
\right) \Bigg|_{y=Y(x)}\\
&\hspace{2cm}
+ Y'(x) \frac{\partial}{\partial y}\left(
(x+y) \ln(x+y) - x \ln(x) - y \ln(y)
\right) \Bigg|_{y=Y(x)}\, .
\end{split}
\end{equation*}
But the first term is $I'(x)$ since we may write
\begin{equation*}
\begin{split}
&\hspace{-1cm}
 \frac{\partial}{\partial x}\left(
(x+y) \ln(x+y) - x \ln(x) - y \ln(y)
\right) \Bigg|_{y=Y(x)}\\
&\hspace{1cm}
=\, 
\left( \frac{\partial}{\partial x_1}\left(
(x_1+y) \ln(x_1+y) - x_1 \ln(x_1) - y \ln(y)
\right) \Bigg|_{x_1=K(y)}\right) \Bigg|_{y=Y(x)}\, .
\end{split}
\end{equation*}
So we have
\begin{equation}
\label{eq:FH}
\begin{split}
&\hspace{-1cm}
\frac{d}{dx}\, \left(
(x+Y(x)) \ln(x+Y(x)) - x \ln(x) 
- Y(x) \ln(Y(x))
\right)\\
&\hspace{1cm}
=\, I'(x)+
Y'(x) \frac{\partial}{\partial y}\left(
(x+y) \ln(x+y) - x \ln(x) - y \ln(y)
\right) \Bigg|_{y=Y(x)}\, .
\end{split}
\end{equation}
But we also have
\begin{equation}
\label{eq:Iformula}
\begin{split}
&\hspace{-1cm}
-I(x) + \left(
(x+Y(x)) \ln(x+Y(x)) - x \ln(x) - Y(x) \ln(Y(x))
\right)\\
&\hspace{1cm}
=\, \lim_{\substack{n \to \infty\\ k/n \to Y(x)}} \frac{1}{n}\, 
\left(\ln\left(\mathcal{M}_n^{(k)}(\boldsymbol{a})\right)- \ln(a_n)\right)\, .
\end{split}
\end{equation}
So, taking derivatives of both sides, and using equation (\ref{eq:FH}), we have
$$
Y'(x) \frac{\partial}{\partial y}\left(
(x+y) \ln(x+y) - x \ln(x) - y \ln(y)
\right) \Bigg|_{y=Y(x)}\, =\,
\frac{d}{dx}
 \lim_{\substack{n \to \infty\\ k/n \to Y(x)}} \frac{1}{n}\, 
\left(\ln\left(\mathcal{M}_n^{(k)}(\boldsymbol{a})\right)- \ln(a_n)\right)\, .
$$
But the right hand side may also be written as 
$$
Y'(x)\, \left(\frac{d}{dy}
 \lim_{\substack{n \to \infty\\ k/n \to y}} \frac{1}{n}\, 
\left(\ln\left(\mathcal{M}_n^{(k)}(\boldsymbol{a})\right)- \ln(a_n)\right)\right)
\Bigg|_{y=Y(x)}\, .
$$
So either $Y'(x)=0$ or else
$$
\frac{\partial}{\partial y}\left(
(x+y) \ln(x+y) - x \ln(x) - y \ln(y)
\right)\, =\, \frac{d}{dy}
 \lim_{\substack{n \to \infty\\ k/n \to y}} \frac{1}{n}\, 
\left(\ln\left(\mathcal{M}_n^{(k)}(\boldsymbol{a})\right)- \ln(a_n)\right)\, ,
$$
where we allow ourselves to drop the prescription $y=Y(x)$ at the end, for ease
of notation.
We will use the latter and we note that it is not a differential equation for $Y(x)$,
but a usual algebraic equation.
(This is the somewhat surprising element if one would have not taken account of the 
fact arising from Fermat's lemma that a certain derivative vanishes.)

Now using the formulas
$$
 \lim_{\substack{n \to \infty\\ k/n \to y}} \frac{1}{n}\, 
\left(\ln\left(\mathcal{M}_n^{(k)}(\boldsymbol{B})\right)- \ln(B_n)\right)\, 
=\, \left(1+\frac{1}{2}\, y\right)\ln\left(1+\frac{1}{2}\, y\right)
-\left(\frac{1}{2}\, y\right)\ln\left(\frac{1}{2}\, y\right)\, ,
$$
and
$$
 \lim_{\substack{n \to \infty\\ k/n \to y}} \frac{1}{n}\, 
\left(\ln\left(\mathcal{M}_n^{(k)}(\boldsymbol{C})\right)- \ln(C_n)\right)\, 
=\, \left(2+ y\right)\ln\left(2+y\right)
-\left(1+y\right)\ln\left(1+y\right) - 2\ln(2)\, ,
$$
coming from Propositions \ref{prop:Catalan} and \ref{prop:central},
the algebraic equations are easily solved to obtain $Y(x)$ (for example, using Wolfram). Then the formulas for $I(x)$ may be obtained from equation (\ref{eq:Iformula}).
This is how we obtained the formulas that we used in Theorem \ref{thm:LDP}.

\section{Proof of Proposition \ref{prop:2SRWBrate} using the generating function}

Let us assume at first that we are in a general set-up with a
$$
g(t)\, =\, \sum_{n=0}^{\infty} a_n t^n\, ,
$$
which is known to us. Let us further assume that
$$
	\lim_{\substack{n \to \infty\\k/n \to \kappa}} \frac{1}{n}\, \ln\left(\mathcal{M}_n^{(k)}(\boldsymbol{a})\right)\, 
=\, \max_{t \in [0,\rho]} \kappa \ln(g(t)) - \ln(t)\, ,
$$
as is normally the case.
This is precisely what we proved in our setting in Theorem \ref{thm:circle}, using a simplified version of the circle method.

Here we assume $\lim_{n \to \infty} \frac{1}{n} \ln(a_n) = \alpha 
\in (-\infty,\infty)$, and  $\rho = e^{-\alpha} \in (0,\infty)$ is the radius of convergence of $g$.
Then for the optimal $t$, namely $T(\kappa)$, we should have
$$
	\kappa\, =\, \frac{g(t)}{t g'(t)}\, \Bigg|_{t=T(\kappa)}\, .
$$
Let us define the inverse of $T$ as $U$, so that
$$
	U(t)\, =\, \frac{g(t)}{t g'(t)}\, .
$$
So we expect to have
$$
\lim_{\substack{n \to \infty\\k/n \to \kappa}} \frac{1}{n}\, \ln\left(\mathcal{M}_n^{(k)}(\boldsymbol{a})\right)\, 
=\, \kappa \ln(g(t)) - \ln(t) \Bigg|_{t=T(\kappa)}\, .
$$
Or in other words, we expect to have
$$
\lim_{\substack{n \to \infty\\k/n \to U(t)}} \frac{1}{n}\, \ln\left(\mathcal{M}_n^{(k)}(\boldsymbol{a})\right)\, 
=\, U(t) \ln(g(t)) - \ln(t)\, .
$$
We call the right-hand-side $W(t)$. 
Using this we may determine that
$$
x(t)\, =\, U(t) \Big(e^{W'(t)/U'(t)}-1\Big)\, .
$$
But
$$
W'(t)\, =\, U'(t) \ln(g(t))\, .
$$
So
$$
x(t)\, =\, U(t)\left(g(t)-1\right)\, .
$$
Using equation (\ref{eq:Iformula}), we may determine that, letting
$$
\widetilde{I}(t)\, =\, I(x(t))\, ,
$$
we have
$$
\widetilde{I}(t)\, =\, -x(t) \ln\left(1-\frac{1}{g(t)}\right)+\ln(t)+\alpha\, .
$$

\section*{Acknowledgments}
We took inspiration from an earlier collaboration with Samen Hossein,
and S.S.~gratefully acknowledges him for useful conversations.
S.S.~is also grateful to Malek Abdesselam for helpful advice, especially for suggesting \cite{CantonFFM}.

\baselineskip=12pt
\bibliographystyle{plain}

\end{document}